\documentclass[12pt,a4paper]{amsart}

\usepackage{amssymb,verbatim}
\usepackage{graphicx,amssymb,amsfonts,epsfig,amsthm,a4,amsmath,url}

\newtheorem{thm}{Theorem}[section]
\newtheorem{cor}[thm]{Corollary}
\newtheorem{lem}[thm]{Lemma}

\newtheorem{prop}[thm]{Proposition}
\theoremstyle{definition}
\newtheorem{defn}[thm]{Definition}
\newtheorem{cla}{Claim}

\newtheorem{rem}[thm]{Remark}
\newtheorem{ex}[thm]{Example}
\numberwithin{equation}{section}

\newcommand{\vv}{\mathfrak{v}}
\newcommand{\mkr}{\mathfrak{r}}
\newcommand{\g}{\mathfrak{g}}
\newcommand{\eps}{\varepsilon}
\newcommand{\Cone}{\textnormal{Cone}}

\begin{document}

\title[Lie groups and cone equivalences]{Asymptotic cones of Lie groups and cone equivalences}\author{Yves de Cornulier}
\address{IRMAR \\ Campus de Beaulieu \\
35042 Rennes Cedex, France}
\email{yves.decornulier@univ-rennes1.fr}
\thanks{Supported by ANR project {\it ``QuantiT"} JC08\textunderscore 318197.}
\date{February 8, 2011}
\maketitle

\begin{abstract}
We introduce cone bilipschitz equivalences between metric spaces. These are maps, more general than quasi-isometries, that induce a bilipschitz homeomorphism between asymptotic cones. Non-trivial examples appear in the context of Lie groups, and we thus prove that the study of asymptotic cones of connected Lie groups can be reduced to that of solvable Lie groups of a special form. We also focus on asymptotic cones of nilpotent groups.
\end{abstract}


\section{Introduction}

Let $X=(X,d)$ be a metric space. The idea of defining an ``asymptotic cone" for $X$, namely a limit when $t\to+\infty$ for the family metric spaces $\frac1tX=(X,\frac1t d)$, was brought in by Gromov \cite{Gro81} in terms of Gromov-Hausdorff convergence. This definition was satisfactory for the purpose of groups with polynomial growth \cite{Pan,Bre}, but to generalize the definition to arbitrary metric spaces, it was necessary to drop the hope of getting a limit in a reasonable topological sense, and consider ultralimits, which make use of the choice of an ultrafilter.
Following Van der Dries and Wilkie \cite{DW}, the limit $\lim_\omega \frac1tX$, formally defined in Section \ref{conemaps}, can be defined when $\omega$ is an unbounded ultrafilter on the set of positive real numbers (by unbounded, we mean that $\omega$ does not contain any bounded subset). This is a metric space, called the asymptotic cone of $X$ with respect to $\omega$, and denoted by $\Cone_\omega(X,d)$ or $\Cone_\omega(X)$ for short. 

Considerable progress in the study of asymptotic cones of groups was then made in Gromov's seminal book \cite{Gro93}, in which the second chapter is entirely devoted to asymptotic cones. Since then, a vast literature appeared on the subject, including the papers \cite{KL,Br,TV,Dru,Ri,KSTT,DrS,Cor,BM}.

The classification of groups in terms of their asymptotic cones is not as fine as the quasi-isometry classification, but in some cases, for instance that of connected Lie groups, it looks more approachable. Here are a few facts relevant to the general study of connected Lie groups up to quasi-isometry.

Let $(\mathcal{C}_0)$ be the class of {\it triangulable} Lie groups, i.e.~groups isomorphic to a closed connected group of real upper triangular matrices. 

\begin{itemize}
\item Every connected Lie group is quasi-isometric to a group in the class $(\mathcal{C}_0)$ \cite[Lemma 6.7]{Cor};
\item open question: is it true that any two quasi-isometric groups in the class $(\mathcal{C}_0)$, are isomorphic? In the nilpotent case, this is considered as a major open question in the field.
\end{itemize}

In the study of the large-scale geometry of a group $G$ in the class $(\mathcal{C}_0)$, a fundamental role is played by the {\it exponential radical} $R$ of $G$, which is defined by saying that $G/R$ is the largest nilpotent quotient of $G$. 

Define $(\mathcal{C}_1)$ as the class of groups $G$ in $(\mathcal{C}_0)$ having a closed subgroup $H$ (necessarily nilpotent) such that
\begin{enumerate}
\item $G$ is the semidirect product $R\rtimes H$;
\item\label{cond2} the action of $H$ on the Lie algebra of $R$ is $\mathbf{R}$-diagonalizable (in particular, $[H,H]$ centralizes $R$).
\end{enumerate}

\begin{rem}Let $G$ be the (topological) unit component of $\mathbf{G}(\mathbf{R})$, where $\mathbf{G}$ is an algebraic $\mathbf{R}$-subgroup of the upper triangular matrices and assume that $\mathbf{G}$ has no nontrivial homomorphism to the additive one-dimensional group (so that the unipotent radical coincides with the exponential radical). Then $G$ is in the class $(\mathcal{C}_1)$.
\end{rem}

In general, and even if $G$ is algebraic, the exact sequence $1\to R\to G\to G/R\to 1$ need not be split \cite[Section~4]{Cor} and there is no subgroup $H$ as above. Even when a splitting exists, Condition (\ref{cond2}) is not satisfied in general. So the class $(\mathcal{C}_1)$ appears as a class of groups, much smaller than the whole class $(\mathcal{C}_0)$, but in which the large-scale geometry is more likely to be understood. Our main result is to associate, to every $G$ in the class $(\mathcal{C}_0)$, a ``nicer" group $G'\in (\mathcal{C}_1)$, of the same dimension and with the same exponential radical, such that the asymptotic cones of $G$ and $G'$ are the same. To state a precise result, we introduce the following new concept. 

We define a map between metric spaces $X,Y$ to be a \textit{cone bilipschitz equivalence} if it induces a bilipschitz homeomorphism at the level of asymptotic cones for \textit{all} unbounded ultrafilters, and we say that $X$ and $Y$ are \textit{cone bilipschitz-equivalent} if there exists such a map. (The precise definitions will be provided and developed in Section \ref{conemaps}.)
It is easy and standard that any quasi-isometry between metric spaces is a cone bilipschitz equivalence, but on the other hand there exist cone bilipschitz-equivalent metric spaces that are not quasi-isometric (see the example before Corollary \ref{coneq_nilp}).

\begin{thm}
Let $G$ be any connected Lie group. Then $G$ is cone bilipschitz equivalent to a group $G_1$ in the class $(\mathcal{C}_1)$. More precisely, if $G$ is triangulable with exponential radical $R$, then there is a split exact sequence
$$1 \to R\to G_1\to G/R\to 1,$$
in which $R$ embeds as the exponential radical of $G_1$.
\end{thm}

If $\g$ is a Lie algebra, denote by $(\g^i)$ its descending central series ($\g^1=\g$; $\g^{i+1}=[\g,\g^i]$); we have $$[\g^i,\g^j]\subset\g^{i+j}$$ for all $i,j$, so the bracket induces a bilinear operation $$[\g^i/\g^{i+1},\g^j/\g^{j+1}]\to \g^{i+j}/\g^{i+j+1}.$$
This defines a Lie algebra structure on $$\g_{\textnormal{grad}}=\bigoplus_{i\ge 1}\g^i/\g^{i+1},$$
called the {\it associated graded Lie algebra}. The Lie algebra $\g$ is called {\it gradable} (or, more commonly but slightly ambiguously, {\it graded}) if it is isomorphic to its associated graded Lie algebra. A simply connected nilpotent Lie group is called gradable if its Lie algebra is gradable. Pansu proved in \cite{Pan2} that any two quasi-isometric gradable Lie groups are actually isomorphic.

The description of asymptotic cones of simply connected nilpotent Lie groups, due to Pansu and Breuillard \cite{Pan,Bre}, can be understood with the help of the notion of cone equivalences. We explain this in detail in Appendix \ref{ap}.
In particular any simply connected nilpotent Lie group $G$ is cone bilipschitz equivalent to its associated graded simply connected nilpotent Lie group $G_{\textnormal{grad}}$ (on the other hand, $G$ and $G_{\textnormal{grad}}$ may be non-quasi-isometric, see Benoist-Shalom's example in \cite[Section~4.1]{Sha04}). This is explicit: if both groups are identified to their Lie algebra through the exponential map, the cone equivalence there is just the identity, the associated graded Lie algebra having the same underlying subspace as the original one. Moreover, for a suitable choice of metrics, the multiplicative constants in the definition of cone equivalence are equal to~1, so they induce isometries between the two asymptotic cones.

\begin{cor}
Any connected Lie group $G$ is cone bilipschitz-equivalent to a group
in the class $(\mathcal{C}_1)$ for which moreover the nilpotent subgroup $H$ (as in the definition of the class $(\mathcal{C}_1)$) is graded.\label{coneq_nilp}
\end{cor}

\medskip

\noindent {\bf Acknowledgments.} I thank Emmanuel Breuillard and Pierre de la Harpe for useful comments.

\setcounter{tocdepth}{1}
\tableofcontents


\section{Cone maps}\label{conemaps}
In this section, $(X,d)$ is a non-empty metric space and $x_0$ is a given point in $X$. We denote, for $x\in X$, $|x|=d(x,x_0)$. The letter $t$ always denotes a variable real number, say $\ge 1$, and typically tending to $+\infty$ (e.g.\ $f(t)\ll t$ means $f(t)/t\to 0$ when $t\to\infty$).

\subsection{Asymptotic cones}\label{asco}
By {\it unbounded ultrafilter} we mean an ultrafilter on the set of positive real numbers, supported at $+\infty$, i.e.\ not containing any bounded subset.

The {\it asymptotic cone} $\Cone_\omega(X,d)$ is defined as the ultralimit with respect to the unbounded ultrafilter $\omega$ of the family of pointed metric spaces $$\left(X,\frac1td,x_0\right)$$ for some $x_0\in X$. This means that 
$\Cone_\omega(X,d)$ is defined as follows: define $\text{Precone}(X,d)$ as the set of families $(x_t)_{t\ge 1}$ in $X$ such that $(d(x_t,x_0)/t)$ is bounded (this set does not depend on $x_0$). Endow it with the pseudo-distance $$d_\omega((x_t),(y_t))=\lim_\omega\frac1td(x_t,y_t),$$
and obtain $\Cone_\omega(X,d)$ as the metric space obtained from the pseudo-metric space
$$(\text{Precone}(X,d),d_\omega)$$ by identifying points at distance zero. After these identifications, the constant sequence $(x)$ does not depend on $x\in X$ and appears as a natural base point for $\Cone_\omega(X,d)$. See Drutu's article \cite{Dru} for a more detailed construction and a survey on asymptotic cones.

\begin{ex}
The metric space $\mathbf{R}$ with the Euclidean distance, is isometric to any of its asymptotic cones, through the map $x\mapsto (tx)$, the reciprocal map being given by $(x_t)\mapsto\lim_\omega\frac{x_t}{t}$.
\end{ex}

Taking asymptotic cones with respect to a given unbounded ultrafilter defines a functor from the category of non-empty metric spaces with large-scale Lipschitz maps to the category of pointed metric spaces with Lipschitz maps preserving base points.

In the following, we extend this functor to a broader class of maps.

\subsection{Cone defined maps}\label{cdm}

For real-valued functions, write $u\preceq v$ if $u\le Av+B$ for some constants $A,B>0$. When applicable, $u\ll v$ means $u(g)/v(g)\to 0$ when $g\to\infty$.

\begin{defn}
Let $Y$ be another non-empty metric space, and also denote by $|\cdot|$ the distance to the given point of $Y$. A map $f:X\to Y$
is {\it cone defined} if for every family $(x_t)$ in $X$, $|x_t|\preceq t$ implies $|f(x_t)|\preceq t$, and moreover for any families $(x_t),(x'_t)$ in $X$ with $|x_t|,|x'_t|\preceq t$ and $d(x_t,x'_t)\ll t$, we have $d(f(x_t),f(x'_t))\ll t$.
\end{defn}

(Note that this does not depend on the choices of the given points.) This is exactly the condition we need to define, for every unbounded ultrafilter $\omega$, the induced map $\Cone_\omega(X)\to\Cone_\omega(Y)$, by mapping the class of a sequence $(x_n)$ to the class of the sequence $(f(x_n))$. The latter map preserves base points of the cones.

\begin{ex}
Let $f:\mathbf{R}\to\mathbf{R}$ be any function satisfying $$\lim_{|x|\to\infty}f(x)/x=0.$$ Then the map $x\mapsto x+f(x)$ is cone defined and induces the identity map at the level of all asymptotic cones.
\end{ex}

\begin{prop}
Let $X,Y$ be metric spaces and $f:X\to Y$ a map. Then $f$ is cone defined if and only if the following two conditions are satisfied
\begin{itemize}
\item $|f(x)|\preceq |x|$;
\item $(d(x,y)+1)/(|x|+|y|)\to 0$ implies $d(f(x),f(y))/(|x|+|y|)\to 0$.
\end{itemize}\label{prop:cone_defined}
\end{prop}
\begin{proof}
The conditions are clearly sufficient. Conversely, if the first one fails, for some sequence $(x_i)$ in $X$, we have $|x_i|\to\infty$ and $|f(x_i)|/|x_i|\to\infty$. 
If $i(t)$ is the largest $i$ such that $|x_i|\le t$, then the sequence $(x_{i(t)})$ satisfies $x_{i(t)}\preceq t$ and $f(x_{i(t)})\npreceq t$. The other part is similar. 
\end{proof}

\begin{prop}
Let $X,Y$ be metric spaces and $f:X\to Y$ a cone defined map. Then, for every unbounded ultrafilter $\omega$, the induced map $\tilde{f}:\Cone_\omega(X)\to\Cone_\omega(Y)$ is continuous.
\end{prop}
\begin{proof}If $x_0$ is the base point of the cone, it follows from the first condition of Proposition \ref{prop:cone_defined} that $d(\tilde{f}(x),\tilde{f}(x_0))\le Cd(x,x_0)$ for some constant $C>0$, for all $x\in X$, so $\tilde{f}$ is continuous at $x_0$.

Suppose $x\in\Cone_\omega(X)-\{x_0\}$ and let us check that $\tilde{f}$ is continuous at $x$. Write $x=(x_t)$ and let $y=(y_t)$. For some function $u$ tending to 0 at 0, we have $d(f(x_t),f(y_t))/(|x_t|+|y_t|)\le u((d(x_t,y_t)+1)/(|x_t|+|y_t|))$ for all $n$. We can suppose that $u\le C$ and is continuous. We have 
$$ \frac{d(f(x_t),f(y_t))}{t}\le \frac{|x_t|+|y_t|}{t}u\left(\frac{d(x_t,y_t)+1}{|x_t|+|y_t|}\right);$$
taking the limit with respect to $\omega$, we obtain
$$d(\tilde{f}(x),\tilde{f}(y))\le (|x|+|y|)u\left(\frac{d(x,y)}{|x|+|y|}\right),$$
which tends to zero when $y$ tends to $x$.
\end{proof}

\begin{defn}Two cone defined maps $f_1,f_2:X\to Y$ are {\it cone equivalent} if the induced maps $\Cone_\omega(X)\to\Cone_\omega(Y)$ are equal for all $\omega$. If $f_2$ is a constant map, we then say that $f_1$ is {\it cone null}.\end{defn}

\begin{prop}
Let $f_1,f_2$ be cone defined maps $X\to Y$. Then $f_1$ and $f_2$ are cone equivalent if and only if for some given point $x_0\in X$ and for some sublinear function $q$, we have, for all $x\in X$
$$d(f_1(x),f_2(x))\le q(|x|).$$ 
\end{prop}
\begin{proof}
The proof is very similar to the previous one. Suppose the condition is satisfied. Let $(x_t)$ be a linearly bounded family in $X$. Then $d(f_1(x_t),f_2(x_t))/t\le q(|x_t|))/t\to 0$ when $t\to +\infty$, so $(f_1(x_t))$ and $(f_2(x_t))$ coincide in any asymptotic cone of $Y$.

Conversely, if the condition is not satisfied, there exists a sequence $(\xi_i)$ tending to infinity in $X$ and $\eps>0$ such that $$d(f_1(\xi_i),f_2(\xi_i))\ge \eps |\xi_i|,\quad\forall i.$$ For every $t$, define $x_t=\xi_i$, where $i$ is chosen maximal so that $|\xi_i|< t+1$ (this is valid as $(|\xi_i|)_i$ goes to infinity). Let $I$ be the set of $t$ for which $x_t\neq x_{t-1}$. This condition means that $x_t=\xi_i$, where $t\le |\xi_i|<n+1$; as $(|\xi_i|)_i$ is unbounded, $I$ is unbounded. Let $\omega$ be a unbounded ultrafilter containing $I$. Then for $t\in I$, we have $$d(f_1(x_t),f_2(x_t))\ge \eps |x_t|\ge \eps t,$$ so $\lim_\omega d(f_1(x_t),f_2(x_t))/t>0$.
\end{proof}

\subsection{Cone Lipschitz maps}

\begin{defn}
A cone defined map $f:X\to Y$ is called a {\em cone $C$-Lipschitz} map if the induced map $\Cone_\omega(X)\to\Cone_\omega(Y)$ is $C$-Lipschitz for all unbounded ultrafilters $\omega$.
\end{defn}

In particular, such a map naturally induces a $C$-Lipschitz map $\Cone_\omega(X)\to\Cone_\omega(Y)$. Say that $f$ is a {\em cone Lipschitz} map if it is cone $C$-Lipschitz for some $C\in [0,+\infty[$. Thus $\Cone_\omega$ is a functor from the category of non-empty metric spaces with cone Lipschitz maps (modulo cone equivalence) to the category of pointed metric spaces with Lipschitz maps preserving base points.

\begin{prop}\label{conelip}
Let $f$ be a map from $X$ to $Y$. Then $f$ is a cone $C$-Lipschitz map if and only if for some sublinear function $q$, we have, for all $x,y\in X$
$$d(f(x),f(y))\le Cd(x,y)+q(|x|+|y|).$$
\end{prop}
\begin{proof}
Suppose that the condition is satisfied. Let $(x_t)$, $(y_t)$ be linearly bounded sequences in $X$. Then $$d(f(x_t),f(y_t))/t-Cd(x_t,y_t)/t\le q(|x_t|+|y_t|))/t,$$ which tends to zero. So $\lim_\omega d(f(x_t),f(y_t))/t-Cd(x_t,y_t)/t\le 0$ for all $\omega$.

Conversely, if the condition is not satisfied, there exists sequences $(\xi_i)$, $(\upsilon_i)$ with $(|\xi_i|+|\upsilon_i|)$ tending to infinity, and $\eps>0$ such that
$$d(f(\xi_i),f(\upsilon_i))-Cd(\xi_i,\upsilon_i)\ge\eps(|\xi_i|+|\upsilon_i|).$$
For every $t$, define $x_t=\xi_i$ and $y_t=\upsilon_i$, where $i=i(t)$ is chosen maximal so that $$|\xi_i|+|\upsilon_i|<t+1$$ (this is valid as $(|\xi_|+|\upsilon_i|)_i$ tends to infinity). It follows that $(x_t)$ and $(y_t)$ are linearly bounded. Let $I$ be the set of $t$ for which $i(t)\neq i(t-1)$. This condition means that for some $i$ (which can be chosen as $i=i(t)$), we have $t\le |\xi_i|+|\upsilon_i|<t+1$; as the sequence $(|\xi_i|+|\upsilon_i|)_i$ is unbounded, $I$ is unbounded. Let $\omega$ be a unbounded ultrafilter containing $I$. Then for $t\in I$, setting $i=i(t)$, we have \begin{align*}d(f(x_t),f(y_t))-Cd(x_t,y_t)= & \; d(f(\xi_i),f(\upsilon_i))-Cd(\xi_i,\upsilon_i)\\
 \ge &\; \eps(|\xi_i|+|\upsilon_i|)\ge\eps t.\qedhere\end{align*}
\end{proof}

\begin{ex}
If $G$ is the Heisenberg group, then the group law: $G\times G\to G$ is cone defined but not cone Lipschitz, see Remark~\ref{hcl}. Looking at the argument, we see that more precisely the map
\begin{eqnarray*}
\mathbf{R}^2 & \to & G  \\
(x,y) &\mapsto & \begin{pmatrix}
  1 & x & 0 \\
  0 & 1 & y \\
  0 & 0 & 1
\end{pmatrix}  =\begin{pmatrix}
  1 & 0 & 0 \\
  0 & 1 & y \\
  0 & 0 & 1
\end{pmatrix}\begin{pmatrix}
  1 & x & 0 \\
  0 & 1 & 0 \\
  0 & 0 & 1
\end{pmatrix}
\end{eqnarray*}
is cone defined but not cone Lipschitz.
\end{ex}

\subsection{Cone bilipschitz maps}

\begin{defn}A cone defined map $f:X\to Y$ is {\it cone} $M$-{\it expansive} if, for every unbounded ultrafilter $\omega$, the induced map $\tilde{f}:\Cone_\omega(X)\to\Cone_\omega(Y)$ is $M$-expansive, i.e. satisfies
$$d(\tilde{f}(x),\tilde{f}(y))\ge Md(x,y)$$
for all $x,y$. The map $f$ is {\it cone expansive} if it is $M$-cone expansive for some $M>0$. The map $f$ is {\it cone} $(C,M)$-{\it bilipschitz} if it is cone $C$-Lipschitz and cone $M$-expansive, and {\it cone bilipschitz} if this holds for some positive reals $M,C$. The map $f$ is {\it cone surjective} if $\tilde{f}$ is surjective for every unbounded ultrafilter, and is called a {\it cone bilipschitz equivalence} if it is both cone surjective and cone bilipschitz.

\end{defn}

The following proposition is proved in the same lines as the previous ones.

\begin{prop}
Let $f$ be a cone map from $X$ to $Y$. Then $f$ is a cone $M$-expansive map if and only if for some sublinear function $\kappa$, we have, for all $x,y\in X$
$$d(f(x),f(y))\ge Md(x,y)-\kappa(|x|+|y|)).\qed$$
\end{prop}

\begin{prop}
Let $f$ be a cone bilipschitz map from $X$ to $Y$ and $y_0$ a base point in $Y$. We have the equivalences
\begin{itemize}
\item[(i)] $f$ is cone surjective;
\item[(ii)] for some sublinear function $c$, we have, for all $y\in Y$, the inequality $d(y,f(X))\le c(|y|)$; 
\item[(iii)] $f$ is a cone bilipschitz equivalence.
\end{itemize}
\end{prop}
\begin{proof}
Clearly (iii) implies (i).

Suppose (ii) and let us prove (iii). Set $c'=c+1$. For any $y\in Y$, choose $x=g(y)$ with $d(f(x),y)\le c'(|y|)$. We claim that $g$ is a cone map. Indeed, $$Md(g(y_1),g(y_2))\le d(f\circ g(y_1),f\circ g(y_2))+\kappa(|g(y_1)|+|g(y_2)|)$$ $$\le d(y_1,y_2)+c'(|y_1|)+c'(|y_2|)+\kappa(|g(y_1)|+|g(y_2)|).$$ 
For some constant $m$, $\kappa(t)\le Mt/2$ for all $t$. Specifying to $y_1=y$ and $y_2=y_0$, we get, for all $y\in Y$
$$M|g(y)|\le M|g(y_0)|+|y|+c'(|y|)+c'(0)+M|g(y)|/2+M|g(y_0)|/2,$$
which gives a linear control of $|g(y)|$ by $|y|$. Thus for some suitable sublinear function $\kappa'$, we have, for all $y_1,y_2$
$$Md(g(y_1),g(y_2))\le d(f\circ g(y_1),f\circ g(y_2))+\kappa'(|y_1|+|y_2|),$$
so that $g$ is cone Lipschitz.

By construction, $f\circ g$ and $\textnormal{Id}_Y$ are cone equivalent. It follows that $f\circ g\circ f$ and $f$ are cone equivalent. Using that $f$ is cone bilipschitz, it follows that $g\circ f$ and $\textnormal{Id}_X$ are cone equivalent.

Finally suppose (ii) does not hold and let us prove the negation of (i). So there exist $\upsilon_i$ in $Y$ and $\eps>0$ with $d(\upsilon_i,f(X))\ge\eps|y_i|$. We can find, by the usual argument, a sequence $(y_n)$ in $Y$ with $|y_n|\le n+1$ for all $n$ and an infinite subset $I$ of integers so that for $n\in I$, $y_n=\upsilon_i$ for some $i$ and $|\upsilon_i|\ge n$. Pick a unbounded ultrafilter on the integers containing $I$. Suppose by contradiction that the sequence $(y_t)$ is image of a sequence $(x_t)$ for the induced map $\Cone_\omega(X)\to\Cone_\omega(Y)$. Then $\lim_\omega d(f(x_t),y_t)/t=0$. But for $t\in I$, we have $d(f(x_t),y_t)=d(f(x_t),\upsilon_i)\ge\eps t$ and we get a contradiction.
\end{proof}

\begin{rem}
The cone Lipschitz (and even large scale Lipschitz map) $f:\mathbf{R}\to\mathbf{R}$ mapping $x$ to $x^{1/3}$ is surjective, however the induced map on the cones is constant, so $f$ is not cone surjective. 
\end{rem}


\section{Illustration: the law in a metric group}
 
Let $G$ be a metric group, i.e.\ endowed with a left-invariant pseudodistance~$d$. As usual, we write $|g|=d(1,g)$. The precone $\textnormal{Precone}(G)$ (see Paragraph \ref{asco}) carries a natural group law, given by $(g_n)(h_n)=(g_n h_n)$. The {\it left} multiplications are isometries of the pseudometric space $\text{Precone}(G)$.

Fix an unbounded ultrafilter $\omega$, and let $\textnormal{Sublin}_\omega(G)$ (resp.\ $\textnormal{Sublin}(G)$) denote the set of families $(g_t)\in\textnormal{Precone}(G)$ such that $|g_t|/t\to 0$ with respect to $\omega$ (resp.\ $|g_t|/t\to 0$ when $t\to+\infty$). This is obviously a subgroup of $\textnormal{Precone}(G)$. Then by definition
$$\Cone_\omega(G)=\textnormal{Precone}(G)/\textnormal{Sublin}_\omega(G).$$
In particular, $\Cone_\omega(G)$ is naturally homogeneous under the action of $\textnormal{Precone}(G)$ by isometries.
In general, $\textnormal{Sublin}_\omega(G)$ is not normal in $\textnormal{Precone}(G)$, and thus $\Cone_\omega(G)$ has no natural group structure. The following proposition gives a simple criterion. Endow $G\times G$ with the $\ell^\infty$ product metric (or any equivalent metric). 

\begin{prop}
Let $G$ be a metric group. We have the equivalences
\begin{itemize}
\item[(i)] the group law $\eta:G\times G\to G$ is cone defined;
\item[(ii)] the group inverse map $\tau:G\to G$ is cone defined;
\item[(iii)] $\textnormal{Sublin}_\omega(G)$ is normal in $\textnormal{Precone}(G)$ for every $\omega$;
\item[(iii')] $\textnormal{Sublin}(G)$ is normal in $\textnormal{Precone}(G)$;
\item[(iv)] for every $e>0$ there exists $\eps(e)>0$ such that for all $g,h\in G$ such that $|h|\ge e|g|$ we have $|g^{-1}hg|+1\ge \eps(e)|g|$. 
\end{itemize}\label{lawcd}
\end{prop}

The proof will be given at the end of this section.

\begin{lem}\label{lcdl}
Let $G$ be a metric group with cone defined law, and $H$ is another metric group with a homomorphic quasi-isometric embedding $H\to G$. Then $H$ also has a cone defined law. In particular, this applies if $G$ is endowed with a word metric with respect to a compact generating subset, and $H\to G$ is a proper homomorphism with cocompact image.
\end{lem}
\begin{proof}This immediate, e.g.\ from Criterion (iii') of Proposition \ref{lawcd}.
\end{proof}

\begin{lem}\label{lcdl2}
If $G$ is a locally compact group with a word metric with respect to a compact generating subset and $H$ is a quotient of $G$, if $G$ has cone defined law, then so does $H$.
\end{lem}
\begin{proof}
We also use Criterion (iii'). We can suppose for convenience that $G$ is endowed with the word length with respect to a compact generating set $S$, and $H$ is endowed with the word length with respect to the image of $S$.
Suppose that $(g_t)\in\textnormal{Precone}(H)$ and $(h_t)\in\textnormal{Sublin}(H)$. We can lift then to elements $\tilde{g}_t$ and $\tilde{h}_t$ of $G$ with $|\tilde{g}_t|=|g_t|$ and $|\tilde{h}_t|=|h_t|$. So by (iii') of Proposition \ref{lawcd}, $|\tilde{g}_t\tilde{h}_t\tilde{g}_t^{-1}|\ll t$, and therefore $|g_th_tg_t^{-1}|\ll t$. 
\end{proof}

\medskip

\begin{ex}Here are some examples of metric groups for which the law is cone defined.
\begin{itemize}
\item[(1)] Arbitrary abelian metric groups are trivial examples, since $\textnormal{Precone}(G)$ is then abelian as well and being normal becomes an empty condition in Criterion (iii') of Proposition \ref{lawcd}.

\item[(2)]Less trivial examples are nilpotent locally compact groups with the word metric with respect to a compact generating subset, see Corollary~\ref{nilc}.

\item[(3)]In particular, using Lemma \ref{lcdl}, it follows from Corollary \ref{nilc} that if a finitely generated group, endowed with a word metric, is finite-by-nilpotent (i.e.\ has a finite normal subgroup with nilpotent quotient), then it satisfies the conditions of Proposition \ref{lawcd}. I~do not know if there are any other examples among finitely generated groups (with a word metric).
\end{itemize}\label{exlacode}
\end{ex}

\begin{ex}\label{ccn}
``Most" metric groups fail to have a cone defined law, as illustrated by the following. In each of these examples, $(x_n)\in\textnormal{Precone}(G)$ and $y_0$ a constant (hence sublinear) sequence, such that $|x_ny_0x_n^{-1}|\simeq n$, so clearly Condition (iii) of Proposition \ref{lawcd} (for instance) fails.
\begin{itemize}
\item $G=\langle a,b\rangle$, the nonabelian free group with the word metric, $x_n=a^n$, $y_0=b$;
\item $G=\langle a,b:a^2=b^2=1\rangle$, the infinite dihedral group with the word metric, $x_n=(ab)^n$, $y_0=a$, $x_ny_0x_n^{-1}=(ab)^{2n}a$; here $G$ has an abelian subgroup of index two, so we see that the converse of Lemma \ref{lcdl} dramatically fails;
\item $G$ the Heisenberg group {\it endowed with the metric induced by inclusion into} $\textnormal{SL}_3(\mathbf{R})$ with its word metric,
$$x_n=\begin{pmatrix}
  1 & 2^n & 0 \\
  0 & 1 & 0 \\
  0 & 0 & 1
\end{pmatrix};\quad y_0=\begin{pmatrix}
  1 & 0 & 0 \\
  0 & 1 & 1 \\
  0 & 0 & 1
\end{pmatrix};\quad x_ny_0x_n^{-1}=\begin{pmatrix}
  1 & 0 & 2^n \\
  0 & 1 & 1 \\
  0 & 0 & 1
\end{pmatrix}.$$ 
In particular we see that Example \ref{exlacode}(2) is specific to the word metric and fails for general nilpotent metric groups.
\end{itemize} 
\end{ex}

\medskip

\begin{proof}[Proof of Proposition \ref{lawcd}]~
\begin{itemize}
\item Suppose (iii) and let us prove (ii). Denote by $\equiv_\omega$ the equality in $\Cone_\omega(G)$. Suppose that $(g_t)\equiv_\omega(g'_t)$, i.e.\ $g'_t=g_ts_t$ with $(s_t)$ $\omega$-sublinear. We have $d(g_t^{-1},{g'_t}^{-1})=|g_ts_t^{-1}g_t^{-1}|$, and $(g_ts_t^{-1}g_t^{-1})$ belongs to $\textnormal{Sublin}_\omega(G)$ since the latter is normal. Thus $(g_t^{-1})\equiv_\omega({g'_t}^{-1})$ and therefore $\tau$ is cone defined.
\item Suppose (ii) and let us prove (i). Suppose that $(g_t)\equiv_\omega(g'_t)$ and $(h_t)\equiv_\omega(h'_t)$. Write $g'_t=g_ts_t$. We have
$$d(g_th_t,g'_th'_t)=|h_t^{-1}s_th'_t|\le |h_t^{-1}s_th_t|+|h_t^{-1}h'_t|.$$
Since $\tau$ is cone defined and $(h_t^{-1}s_t^{-1})\equiv_\omega(h_t^{-1})$, we have $(s_th_t)\equiv_\omega(h_t)$, hence by left multiplication (which is isometric) $(h_t^{-1}s_th_t)\equiv_\omega(1)$, i.e.\ $(h_t^{-1}s_th_t)$ is $\omega$-sublinear. Since this holds for any $\omega$, the group law $\eta:G\times G\to G$ is cone defined;
\item Suppose (i) and let us prove (iii). Suppose that $(g_t)\in\textnormal{Precone}(G)$ and $(s_t)\in\textnormal{Sublin}_\omega(G)$. Then, since $(1)\equiv_\omega(s_t)$ and $\eta$ is cone defined, we have $(\eta(1,g_t))\equiv_\omega(\eta(s_t,g_t))$, that is, $|g_t^{-1}s_tg_t|$ is $\omega$-sublinear. So $\textnormal{Sublin}_\omega(G)$ is normal for any $\omega$.
\item Suppose (iii) and let us check (iii'). If $\textnormal{Sublin}_\omega(G)$ is normal for any $\omega$, then so is $\textnormal{Sublin}(G)=\bigcap_\omega\textnormal{Sublin}_\omega(G)$.

\item Suppose that (iv) fails and let us show that (iii') fails. There exists $e_0>0$ and sequences $(g_n),(h_n)$ in $G$ such that $|h_n|\ge e_0|g_n|$ and $(|g_n^{-1}h_ng_n|+1)/|g_n|\to 0$. Note that this forces $|g_n|\to\infty$. Since $|g_n^{-1}h_ng_n|\ge |h_n|-2|g_n|$, we see that $|h_n|\le 3|g_n|$ for $n\ge n_0$ large enough. 
If $r\ge r_0=|g_{n_0}|$, set $n(r)=\sup\{n\ge n_0:|g_n|\le r\}$, $u_r=g_{n(r)}$ and $v_r=h_{n(r)}$. For all $r\ge r_0$ we have $|u_r|\le r$, $e_0|u_r|\le|v_r|\le 3r$, and the family $(u_r^{-1}v_ru_r)$ is sublinear and both families $(u_r)$ and $(v_r)$ are at most linear. Set $N=\{n:\forall m>n,|g_m|>|g_n|\}$. If $r$ belongs to the unbounded set $K=\{|g_n|:n\in N\}$ we have $|u_r|=r$ for all $r\in K$ and therefore the family $(v_r)$ is not sublinear.

\item Suppose that (iv) holds and let us prove (iii). Fix $\omega$. Suppose that $(g_t)\in\textnormal{Precone}(G)$ and $(s_t)\in\textnormal{Sublin}_\omega(G)$. If $(g_ts_tg_t^{-1})\notin\textnormal{Sublin}_\omega(G)$, then for some $I\in\omega$, some $e_0>0$ and all $t\in I$ we have $|g_ts_tg_t^{-1}|\ge e_0|g_t|$. So $$|s_t|+1=|g_t^{-1}(g_ts_tg_t^{-1})g_t|+1\ge\eps(e_0)|g_t|,\quad\forall t\in I,$$
and thus $|g_t|\in\textnormal{Sublin}_\omega(G)$, and in turn $$(g_ts_tg_t^{-1})\in\textnormal{Sublin}_\omega(G),$$ a contradiction.\qedhere
\end{itemize}
\end{proof}


\section{Cone equivalences between Lie groups}

\subsection{Reduction to the split case with semisimple action}

Let $G$ be a triangulable group. We shamelessly use the identification between the Lie algebra and the Lie group through the exponential map.
Precisely, $G$ has two laws, the group multiplication, and the addition $g+h$, which formally denotes $\exp(\log(x)+\log(y))$, using that the exponential function is duly a homeomorphism for a triangulable group \cite{Dix}.
Let $R$ be the exponential radical of $G$, $H$ a Cartan subgroup \cite[Chap.~7, pp.19-20]{Bk}. This means that $H$ is nilpotent and $RH=G$. Set $W=H\cap R$, and let $V$ be a complement subspace of $W$ in $H$ (viewed as Lie algebras).

We denote the word length (with respect to some compact generating subset) in a group by $|\cdot|$, by default the group is $G$, otherwise we make it explicit, e.g. $|\cdot|_H$ denotes word length in $H$.
If $h\in H$, it can be written in a unique way as $rv$ with $r\in W$ and $v\in V$, we write $r=\delta(h)$ and $v=[h]$.

The following lemma is obtained in the proof of \cite[Theorem~5.1]{Cor}.
\begin{lem}
For $v\in V$ we have 
$$|v|\simeq |v|_{G/R}.\qed$$\label{compnorm}
\end{lem}

\begin{lem}\label{majdelta}
For $v,v'\in V$ we have
$$|\delta(v^{-1}v')|\preceq \log(1+|v|_{G/R})+\log(1+|v'|_{G/R}).$$
\end{lem}
\begin{proof}
As $\delta(v^{-1}v')=v^{-1}v'[v^{-1}v']^{-1}$,
$$|\delta(v^{-1}v')|_H\le |v|_H+|v'|_H+|[v^{-1}v']|_H,$$
so by Lemma \ref{compnorm}
$$|\delta(v^{-1}v')|_H\preceq |v|_{G/R}+|v'|_{G/R}+|v^{-1}v'|_{G/R}\preceq |v|_{G/R}+|v'|_{G/R}.$$
Since $\delta(v^{-1}v')$ belongs to the exponential radical,
$$|\delta(v^{-1}v')|\preceq\log(1+|\delta(v^{-1}v')|_H)$$ and we get the conclusion.
\end{proof}

\begin{lem}For $r\in R,v\in V$,
$$|rv|\simeq |r|+|v|_{G/R}.$$\label{estlength}\end{lem}
\begin{proof}

Clearly $|v|_{G/R}\preceq |v|$, so we obtain the inequality $\preceq$.

From the projection $G\to G/R$, we get $|v|_{G/R}\preceq |rv|$.
On the other hand $|r|\le |rv|+|v|$. Now $|v|\preceq |v|_{G/R}$ by Lemma \ref{compnorm}. So $|rv|\succeq |r|+|v|_{G/R}$.
\end{proof}

It is wary here to distinguish the Lie group and its Lie algebra. Thus denote by $\mkr$ the Lie algebra of $R$ and consider the action of $H$ on $\mkr$, given by restriction of the adjoint representation. It is given by a homomorphism $\alpha$ from $H$ to the automorphism group of the Lie algebra $\mkr$, which is an algebraic group. This action is triangulable, and we can write in a natural way, for every $h\in H$, $\alpha(h)=\beta(h)u(h)$ where $\beta$ is the diagonal part and $u(h)$ is the unipotent part; both are Lie algebra automorphisms of $\mkr$ (view this by taking the Zariski closure of $\alpha(H)$, which can be written in the form $DU$ with $D$ a maximal split torus and $U$ the unipotent radical), and $\beta$ defines a continuous action of $H$ on $\mkr$. Note that this action is trivial on $[H,H]$ and hence on $W$.

As $H$ is nilpotent, we can decompose $\mkr$ into characteristic subspaces for the $\alpha$-action; this way we see in particular that $u$ is an action as well (this strongly relies on the fact that $H$ is nilpotent and connected). In such a decomposition, we see that the matrix entries of $u(h)$, for $h\in H$, are polynomially bounded in terms of $|h|$ (more precisely, $\le C|h|^k$, where $C$ is a constant (depending only on $G$ and the choice once and for all of a basis adapted to the characteristic subspaces) and $k+1$ is the dimension of $\mkr$.

Now define $A(h)$, $B(h)$, resp. $U(h)$, as the automorphism of $R$ whose tangent map is $\alpha(h)$, $\beta(h)$, resp. $u(h)$. Note that $\alpha(h)$ is the left conjugation by $h$ in $R$. Then $B$ provides a new action of $H$ on $R$. We consider the group $R\rtimes H/W$ defined by the action $B$. We write it by abuse of notation $R\rtimes V$, identifying $V$, as a group, to $H/W$.

Consider the map $$\begin{array}{cc}\psi: & G=RV\to R\rtimes V \\ & rv\mapsto rv\end{array}.$$

In general $\psi$ is not a quasi-isometry (and is even not a coarse map, i.e. there exists a sequence of pairs of points at bounded distance mapped to points at distance tending to infinity).

\begin{thm}\label{psi}
The map $\psi$ is a cone bilipschitz equivalence.
\end{thm}

\begin{lem}
If $a,b\ge 0$ and $c\ge 1$, and if $|\log(a)-\log(b)|\le\log(c)$, then $|\log(1+a)-\log(1+b)|\le\log(c)$.\qed\label{abc}
\end{lem}

\begin{proof}[Proof of Theorem \ref{psi}]
On both groups, the word length is equivalent by Lemma \ref{estlength} to $L(rv)=\ell(r)+|v|_{G/R}$, where $\ell$ is a length (subadditive and symmetric) on $R$ with $\ell(r)\simeq\log(1+\|r\|)$, and we consider the corresponding left-invariant ``distances"\footnote{These are not necessarily distances but are equivalent to the left-invariant word or Riemannian distances.} $d_1$ and $d_2$.

We consider two group laws at the same time on the Cartesian product $R\times V$; in order not to introduce tedious notation, we go on writing both group laws by the empty symbol; on the other hand to avoid ambiguity, we write the length $L$ as $L=L_1=L_2$ (and $\ell=\ell_1=\ell_2$) and when we compute a multiplication inside the symbols $L_1()$ or $d_1(,)$, we mean the multiplication inside $G$, while in $L_2$ and $d_2$ we mean the new multiplication from $R\rtimes V$.

We have, for $r,r'\in R;v,v'\in V$
$$d_1(rv,r'v')=L_1(v^{-1}r^{-1}r'v')=L_1(A(v^{-1})(r^{-1}r')v^{-1}v')$$
We have $$d_1(rv,r'v')=L_1(A(v^{-1})(r^{-1}r')\delta(v^{-1}v')[v^{-1}v'])$$ $$=\ell_1(A(v^{-1})(r^{-1}r')\delta(v^{-1}v'))+|v^{-1}v'|_{G/R}$$
and
$$d_2(rv,r'v')=\ell_2(B(v^{-1})(r^{-1}r'))+|v^{-1}v'|_{G/R}$$
(no $\delta$-term for $d_2$ because $V$ is a subgroup for the second law).

So $$d_1(rv,r'v')-d_2(rv,r'v')$$ $$=\ell_1(A(v^{-1})(r^{-1}r')\delta(v^{-1}v'))-\ell_2(B(v^{-1})(r^{-1}r')).$$
Write for short $\rho=B(v^{-1})(r^{-1}r')$ (which defines the same element for the two group laws), so
$$d_1(rv,r'v')-d_2(rv,r'v')=\ell_1(U(v^{-1})(\rho)\delta(v^{-1}v'))-\ell_2(\rho).$$
We then get
$$|d_1(rv,r'v')-d_2(rv,r'v')|=|\ell(U(v^{-1})(\rho))-\ell(\rho)|+\ell(\delta(v^{-1}v')).$$
(We see that we can write $\ell$ because this does no longer depends on the choice of one of the two laws.)

Let us now work on the manifold $R$, endowed with a Riemannian length $\lambda$. As the tangent map of $U(v^{-1})$ is $C|v|^k$-bilipschitz (i.e. both it and its inverse are $C|v|^k$-Lipschitz), the same is true for $U(v^{-1})$ itself. In particular,
$$\frac{1}{C|v|^k}\lambda(\rho)\le\lambda(U(v^{-1})\rho)\le C(1+|v|)^k\lambda(\rho).$$
This means that $$|\log(\lambda(\rho))-\log(U(v^{-1})\rho))|\le \log(C)+k\log(1+|v|).$$
By Lemma \ref{abc} (we can pick $C\ge 1$), $$|\log(1+\lambda(\rho))-\log(1+U(v^{-1})\rho))|\le \log(C)+k\log(1+|v|),$$
thus
$$|\ell(\rho)-\ell(U(v^{-1})\rho))|\preceq \log(1+|v|);$$
moreover by Lemma \ref{majdelta}
$$\delta(v^{-1}v')\preceq \log(1+|v|)+\log(1+|v'|).$$
Accordingly,
$$|d_1(rv,r'v')-d_2(rv,r'v')|\preceq \log(1+|v|)+\log(1+|v'|)\ll |rv|+|r'v'|.$$
Thus the map $\psi$ is a cone equivalence (with constants one). This means that the cones defined from $d_1$ and $d_2$ are isometric; however note that these are equivalent to metrics but are not necessarily metrics (they are maybe not subadditive) so the statement obtained for the usual cones (defined with genuine metrics) are that they are bilipschitz.
\end{proof}


\appendix

\section{Asymptotic cones of nilpotent groups, after Pansu and Breuillard}\label{ap}

Fix a $d$-dimensional real $s$-nilpotent Lie algebra $\g$. It is endowed with a group structure denoted by no sign or by a dot $(\cdot)$, defined by the Baker-Campbell-Hausdorff formula. We sometimes write $G=(\g,\cdot)$, but the underlying sets are the same.

For all $i$, denote by $\vv^i$ a complement subspace of $\g^{i+1}$ in $\g^i$. In particular,
\begin{equation}\g=\g^1=\bigoplus_{i=1}^s \vv^i.\label{decompo}\end{equation}

Define, for any $t\in\mathbf{R}$, the endomorphism $\delta_t$ of the graded vector space $\g=\bigoplus\vv^i$ whose restriction to $\vv^i$ is the scalar multiplication by $t^i$. Observe that $t\mapsto\delta_t$ is a semigroup action of $(\mathbf{R},\cdot)$.

\subsection{The cone in the graded case} 
 
We assume {\it in this paragraph} that $(\vv^i)$ is a grading of the Lie algebra $\g$, i.e.\ $$[\vv^i,\vv^j]\subset \vv^{i+j}\quad\forall i,j.$$
It follows that for all $t$, $\delta_t$ is a Lie algebra endomorphism, and therefore is also an endomorphism of $G$.

Consider the subspace $\vv^1$ of $\g$ (endowed with some norm), viewed as the tangent space of $G$ at $1$, and translate it by left multiplication, providing for any $g\in G$ a subspace $\mathcal{D}_g$ of the tangent space $\mathcal{T}_g$ along with a norm, in a way which is compatible with left multiplication. The Carnot-Carath\'eodory distance between $x$ and $y$ is the following ``sub-Finslerian" distance
$$d_{\textnormal{CC}}(x,y)=\inf_\gamma L(\gamma),$$ where $\gamma$ ranges over regular smooth paths everywhere tangent to $\mathcal{D}$ and $L(\gamma)$ is the length of $\gamma$, computed by integration, using the norm on each $\mathcal{D}_g$. This distance is finite (i.e.\ every two points can be joined by such a path), defines the usual topology of $G$ \cite[Theorem~1]{Be}; moreover it is equivalent to the word length (equivalent here refers to large scale equivalence). 

It is immediate from the definition of $d_{\textnormal{CC}}$ that
$$d_{\textnormal{CC}}(\delta_t(x),\delta_t(y))=td_{\textnormal{CC}}(x,y)\quad \forall x,y\in G;$$
in particular, for any unbounded ultrafilter $\omega$, the map
\begin{eqnarray*}j:(G,d_{\textnormal{CC}})&\to &\Cone_\omega(G,d_{\textnormal{CC}}))\\
g & \mapsto & (\delta_t(g))
\end{eqnarray*}
is an isometric embedding; moreover since balls in $(G,d_{\textnormal{CC}})$ are compact, it follows that this is a surjective isometry, whose inverse is given by $(g_t)\mapsto\lim_\omega\delta_{1/t}(g_t)$. 
This gives in particular the following statement.

\begin{prop}\label{grasy}
Let $G$ be any gradable simply connected nilpotent Lie group. Then
\begin{itemize}
\item[(i)] any asymptotic cone of $G$ is naturally bilipschitz to $(G,d_{\textnormal{CC}})$ (with constants independent of $\omega$);
\item[(ii)] the group law $\eta:G\times G\to G$ is cone defined;
\item[(iii)] under the identification of $G$ with its cone, we have $\tilde{\eta}=\eta$.
\end{itemize}
\end{prop}

\begin{cor}\label{nilc}
Let $G$ be an arbitrary nilpotent compactly generated locally compact group (e.g.\ discrete, or connected). Then the law $G\times G\to G$ is cone defined.
\end{cor}
\begin{proof}
If $G$ is a free $s$-nilpotent group, this follows from Proposition \ref{grasy}. By Lemma \ref{lcdl2}, to have cone defined law is stable under taking quotients, so if $G$ is an arbitrary simply connected nilpotent Lie group then $G$ has cone defined law. 

If $G$ is an arbitrary nilpotent compactly generated locally compact group, then $G$ has a unique maximal compact subgroup $K$; this is a characteristic subgroup and the quotient $H=G/K$ is a torsion-free compactly generated nilpotent Lie group (see \cite[p.~104]{GKR}), i.e.\ $H_0$ is a simply connected nilpotent Lie group and the discrete group $H/H_0$ is a finitely generated torsion-free nilpotent group. By a result essentially due to Malcev (see \cite{Wa} for the general case, when $H$ is not necessarily discrete) $H$ embeds as a closed cocompact subgroup in a simply connected nilpotent Lie group $\tilde{H}$, which has cone defined law by the preceding case.
By Lemma \ref{lcdl}, it follows that $H$ and then $G$ has cone defined law.
\end{proof}

\medskip

\begin{proof}[Proof of Proposition \ref{grasy}]~
\begin{itemize}
\item[(i)] if $|\cdot|$ is a word length on $G$, then the identity map $(G,|\cdot|)\to (G,d_{\textnormal{CC}})$ is a quasi-isometry, so the induced map $$\Cone_\omega(G,|\cdot|)\to\Cone_\omega(G,d_{\textnormal{CC}})=(G,d_{\textnormal{CC}})$$
is bilipschitz, with constants depending only on the quasi-isometry constants, not on $\omega$.

\item[(ii)] We have to check that if $(g_t)=(g'_t)$ and $(h_t)=(h'_t)$ in $\Cone_\omega(G)$, then $(g_th_t)=(g'_th'_t)$ as well. Indeed,
$$\delta_{1/t}((g_th_t)^{-1}g'_th'_t)=\delta_{1/t}(h_t)^{-1}\delta_{1/t}(g_t)^{-1}\delta_{1/t}(g'_t)\delta_{1/t}(h'_t);$$
since $\delta_{1/t}(h_t)^{-1}$ and $\delta_{1/t}(h'_t)$ are bounded (hence convergent for $\omega$), 
\penalty-10000$\delta_{1/t}(g_t)^{-1}\delta_{1/t}(g'_t)$ tends to zero and by compactness, the product tends to zero by continuity of the law in $G$, i.e.\ $\delta_{1/t}((g_th_t)^{-1}g'_th'_t)\to 0$, which means that $$d_{\textnormal{CC}}(1,g_th_t)^{-1}g'_th'_t))=td_{\textnormal{CC}}(1,\delta_{1/t}((g_th_t)^{-1}g'_th'_t))\ll t,$$
using that $d_{\textnormal{CC}}$ defines the usual topology on $G$.
\item[(iii)] Since we know the law is well defined, it is enough to compute in on $G$: the new product of $g$ and $h$ is by definition $$j^{-1}(j(g)j(h))=j^{-1}((\delta_t(g))(\delta_t(h)))$$ $$=j^{-1}((\delta_t(g)\delta_t(h)))=j^{-1}((\delta_t(gh)))=gh.\qedhere$$
\end{itemize}
\end{proof}

\begin{rem}\label{hcl}
If $G$ is a graded simply connected nilpotent Lie group, its group law is not, in general, cone Lipschitz, although it is cone defined by Proposition \ref{grasy}(i). By Proposition \ref{grasy}(iii), this amounts to prove that in $(G,d_{\textnormal{CC}})$, the group law is not Lipschitz. Let us check it when $G$ is the Heisenberg group, the argument probably extending to an arbitrary nonabelian $G$. Set
$$x_n=\begin{pmatrix}
  1 & n & 0 \\
  0 & 1 & 0 \\
  0 & 0 & 1
\end{pmatrix};\quad y_0=\begin{pmatrix}
  1 & 0 & 0 \\
  0 & 1 & 1 \\
  0 & 0 & 1
\end{pmatrix};\quad x^{-1}_ny_0x_n=\begin{pmatrix}
  1 & 0 & -n \\
  0 & 1 & 1 \\
  0 & 0 & 1
\end{pmatrix}.$$ 
Defining $R_x:G\to G$ by $R_x(g)=gx$, we have $$d_{\textnormal{CC}}(R_{x_n}(1),R_{x_n}(y_0))=d_{\textnormal{CC}}(1,x^{-1}_ny_0x_n)\simeq\sqrt{n};$$
so the Lipschitz constant of $R_{x_n}$ is $\succeq\sqrt{n}$ and since each $R_{x_n}$ is the restriction of the group law to a certain subset of $G\times G$, the group law itself cannot be Lipschitz.
\end{rem}

\subsection{The asymptotic cone for general nilpotent Lie groups}\label{dac}

Now we consider the graded Lie algebra associated to the decomposition $\g=\bigoplus \vv^i$. The new bracket is given by
$$[x,y]^\infty=\sum_{i,j}[x_i,y_j]_{i+j},$$
and we denote by $\boxtimes$ the corresponding group law, provided by the Baker-Campbell-Hausdorff formula \cite{Ha}. If $x=\sum_{i\ge 1} x_i$ in this decomposition, define
$$\sigma(x)=\sum_{i\ge 1} \|x_i\|^{1/i}.$$

We call the groups $G=(\g,\cdot)$ and $G_\infty=(\g,\boxtimes)$.

Let $d$ be any left-invariant pseudometric on $G$, quasi-isometric to the word distance with respect to some compact generating subset $S$, i.e.\  $|\cdot|_S\preceq |\cdot|\preceq |\cdot|_S$ (see Paragraph \ref{cdm}).

Given an unbounded ultrafilter $\omega$ on the positive real numbers, for any $g,h\in\g$, define the limit
$$d_\omega(g,h)=\lim_{t\to\omega}\frac{d(\delta_t(g),\delta_t(h))}{t}.$$

\begin{thm}\label{mbp}
The distance $d_\omega$ is a continuous, $\boxtimes$-left-invariant distance on $G$. Moreover, the identity map $\g\to\g$ is a cone 1-bilipschitz equivalence $(G_\infty,d_\omega)\to (G,d)$, thus inducing a bijective isometry 
\begin{eqnarray*}(G_\infty,d_\omega)&\to & \Cone_\omega(G,d)\\
g &\mapsto &(\delta_t(g)).
\end{eqnarray*}
Besides, if $\eta$ is the group law in $G$, then $\eta$ and $\boxtimes$ are cone defined on $(G,d)$ and are cone equivalent; in particular the above isometry $G_\infty\to\Cone_\omega(G)$ is a group isomorphism.
\end{thm}

Before proving Theorem \ref{mbp}, let us give a corollary.
Following Breuillard \cite{Bre}, the metric $d$ is {\it asymptotically geodesic} if for every $\eps$ there exists $s$ such that for all $x,y\in G$ there exist $n$ and $x=x_0,\dots,x_n=y$ such that $d(x_i,x_{i+1})\le s$ for all $i$ and $\sum_i d(x_i,x_{i+1})\le (1+\eps)d(x,y)$. It is straightforward that if the metric $d$ is asymptotically geodesic, then $\Cone_\omega(\g,d)$ is geodesic.

\begin{cor}\label{ccb}
If $d$ is asymptotically geodesic, then $d_\omega$ is the Carnot-Carath\'eodory metric associated to a supplement subspace $\mathfrak{v}_\omega$ of $[\g_\infty,\g_\infty]$ in $\g_\infty$.
\end{cor}

\begin{rem}
If $d$ is asymptotically geodesic, Breuillard \cite{Bre} proves in addition that $d_\omega$ is independent of $\omega$, or equivalently that $\vv_\omega$ is independent of $\omega$.

On the other hand, if $d$ is not asymptotically geodesic, then $d_\omega$ may depend on $\omega$. For instance, in the real Heisenberg group (viewed as a group of upper unipotent 3-matrices), consider the word metric associated with the following weighted generating subset: elements with coefficients of absolute value $\le 1$, with weight 1, and the elements $\begin{pmatrix}1 & 0 & (2n)!^2\\0 & 1 & 0\\ 0 & 0 & 1\end{pmatrix}$ with weight $(2n)!/2$. This length is equivalent to the word metric with respect to a compact generating subset. However, the length of $\begin{pmatrix}1 & 0 & (2n)!^2\\0 & 1 & 0\\ 0 & 0 & 1\end{pmatrix}$ is $(2n)!/2$ while the length of $\begin{pmatrix}1 & 0 & (2n+1)!^2\\0 & 1 & 0\\ 0 & 0 & 1\end{pmatrix}$ is approximately $(2n+1)!/\sqrt{2}$, so $d_\omega\left(\begin{pmatrix}1 & 0 & 1\\0 & 1 & 0\\ 0 & 0 & 1\end{pmatrix},1\right)$ does depend on $\omega$.
\end{rem}

\begin{proof}[Proof of Theorem \ref{mbp}]
Set $\rho(g)=d(0,g)$ and $\rho_\omega(g)=d_\omega(0,g)$. Since $d$ is quasi-isometric to the word distance, it is part of \cite[Proof of Th\'eor\`eme~II.2']{Gu} that we have, for suitable constants and for all $g\in\g$
\begin{equation}\label{cosi}C_1\sigma(g)-C'_1\le \rho(g)\le C_2\sigma(g)+C'_2.\end{equation} 
Observe that $\sigma(\delta_t(g))=|t|\sigma(g)$. Therefore, for $t>0$ 
$$C_1\sigma(g)-\frac{C'_1}{t}\le \frac{\rho(\delta_t(g))}{t}\le C_2\sigma(g)+\frac{C'_2}t;$$ 
in particular we deduce 
\begin{equation}C_1\sigma(g)\le \rho_\omega(g)\le C_2\sigma(g).\label{contrrho}\end{equation} 

\begin{cla}\label{rhorhorho}
For any $(g_t),(h_t)$ in $\g$ ($t\ge 0$) with $\sigma(g_t)+\sigma(h_t)\preceq t$, we have $$\lim_{t\to+\infty}\frac{\rho(g_t^{-1}h_t)-\rho(g_t^{-1}\boxtimes h_t)}{t}= 0.$$
\end{cla}
Indeed, we have 
$$|\rho(g_t^{-1}h_t)-\rho(g_t^{-1}\boxtimes h_t)|\le\rho((g_t)^{-1}h_t)^{-1}(g_t^{-1}\boxtimes h_t))$$
and by Proposition \ref{estim} (and in view of (\ref{cosi})), \begin{equation}\rho((g_t^{-1}h_t)^{-1}(g_t^{-1}\boxtimes h_t))\ll t,\label{rhoeta}\end{equation}
and the claim is proved.

By Claim \ref{rhorhorho}, we have
\begin{align*}d_\omega(g,h) = &  \lim_\omega\frac{\rho(\delta_t(g)^{-1}\delta_t(h))}{t}\\
= & \lim_\omega\frac{\rho(\delta_t(g)^{-1}\boxtimes\delta_t(h))}{t}=\rho_\omega(g^{-1}\boxtimes h),\end{align*}
so $d_\omega$ is $\boxtimes$-left-invariant.
As a pointwise limit of pseudodistances, $d_\omega$ is a pseudodistance. In particular, $\rho_\infty$ is a length on $(\g,\boxtimes)$. Since by (\ref{contrrho}), $\rho_\omega$ is continuous at zero, it immediately follows that $\rho_\omega$ (hence $d_\omega$) is everywhere continuous.

Let us turn to the part concerning asymptotic cones. First observe that the map
\begin{eqnarray}i:(G_\infty,d_\omega)&\to & \Cone_\omega(G_\infty,d_\infty)\label{imap}\\
g &\mapsto &(\delta_t(g))\notag
\end{eqnarray}
is a well-defined isometric embedding, by definition of $d_\omega$.
Moreover, it is surjective, its inverse being given by $j:(g_t)\mapsto\lim_\omega\delta_{1/t}(g_t)$. To check that $i\circ j$ is the identity, compute
\begin{align*}
\tilde{d}_\omega((g_t),(i\circ j(g_t)))=& \lim_{t\to\omega}\frac{1}td_\omega\left(g_t,\delta_t\left(\lim_{u\to \omega}\delta_{1/u}(g_u)\right)\right)\\
=& \lim_{t\to \omega} d_\omega\left(\delta_{1/t}(g_t),\lim_{u\to \omega}\delta_{1/u}(g_u)\right)\\
=& d_\omega\left(\lim_{t\to\omega}\delta_{1/t}(g_t),\lim_{u\to \omega}\delta_{1/u}(g_u)\right)=0,
\end{align*}
where we use continuity of $d_\omega$ in the last line.

We know that $\eta$ is cone-defined by Corollary \ref{nilc}; $\eta$ and $\boxtimes$ are cone-equivalent as a consequence of (\ref{rhoeta}).\end{proof}

\begin{proof}[Proof of Corollary \ref{ccb}]
Since $(G_\infty,d_\omega)$ is geodesic, by \cite[Theorem~2]{Be}, it is the Carnot-Carath\'eodory metric defined by a generating subspace $\vv_\omega$ of the Lie algebra $\g_\infty=(\g,[\cdot,\cdot]^\infty)$. In particular, $\vv_\omega+[\g_\infty,\g_\infty]=\g_\infty$. By (\ref{contrrho}) in the proof of Theorem \ref{mbp}, $([G_\infty,G_\infty],d_\omega)$ cannot contain a bilipschitz copy of a segment, and therefore $\vv_\omega\cap[\g_\infty,\g_\infty]=\{0\}$.
\end{proof}

\subsection{Computation in nilpotent groups} 
In all this paragraph, $\g$ is a finite-dimensional real $s$-nilpotent Lie algebra, on which a norm has been fixed; multiplying the norm by a suitable scalar if necessary we suppose that $\|[x,y]\|\le\|x\|\|y\|$ for all $x,y$.

The group law can be written down by the Baker-Campbell-Hausdorff formula as
\begin{equation}\label{bch}xy=x+y+\sum_{k=2}^s \sum_{u\in\mathcal{U}(k-2)}\lambda_u [u_1,\dots,u_{k-2},x,y],\end{equation}
where $\mathcal{U}(k)=\mathcal{U}(k)[x,y]$ is the set (of cardinality $2^{k-2}$) of functions $\{1,\dots,k-2\}\to\{x,y\}$ ($x,y$ denoting two distinct symbols) and $\lambda_u$ are rational numbers (which can be fixed once and for all, independently of $G$ and $s$).

Besides, for $n\ge 2$ and $x^1,\dots,x^n\in\g$, define by induction $[x^1,\dots,x^n]$ as the usual bracket for $n=2$, and as $[x^1,[x^2,\dots,x^n]]$ for $n\ge 3$. Also by $x^j_i$ we mean the $i$th component of $x^j$ with respect to the linear grading $\g=\bigoplus_i\mathfrak{v}^i$.   
We first need the following computational lemma, which is a formal consequence of multilinearity.
 
\begin{lem}\label{expa}
For any $x^1,\dots,x^k\in\g$, we have
$$[x^1,\dots,x^k]_\ell=\sum_{i_1+\dots +i_k\le \ell}[x^1_{i_1},\dots,x^k_{i_{k}}]_\ell$$
and
$$[x^1,\dots,x^k]^\infty_\ell=\sum_{i_1+\dots +i_k= \ell}[x^1_{i_1},\dots,x^k_{i_{k}}]_\ell.$$
\end{lem}
\begin{proof}
The first equality is obtained by writing $x^j=\sum_{i_j}x^j_{i_j}$ and expanding, noting that $[x^1_{i_1},\dots,x^k_{i_{k}}]_\ell=0$ if $\ell <i_1+\dots +i_k$. For the second, it is enough to observe that $[x^1_{i_1},\dots,x^k_{i_{k}}]^\infty$ is the projection of $[x^1_{i_1},\dots,x^k_{i_{k}}]$ on $\vv^\ell$, for $\ell=i_1+\dots +i_k$.
\end{proof}

\begin{lem}\label{compar}
We have in $\g$, for any integer $\ell\ge 1$
$$\sup_{x,y}\frac{\|(xy-x\boxtimes y)_\ell\|}{\max(\sigma(x),\sigma(y))^{\ell-1}}<\infty.$$
\end{lem}
\begin{proof}
Let us write the Baker-Campbell-Hausdorff formula as in (\ref{bch}).
\begin{align*} & xy-x\boxtimes y\\
& =\sum_{k=2}^s \sum_{u\in\mathcal{U}(k-2)}\lambda_u ([u_1,\dots,u_{k-2},x,y]-[u_1,\dots,u_{k-2},x,y]^\infty),\end{align*}
so, using Lemma \ref{expa},
\begin{align*} & (xy-x\boxtimes y)_\ell\\ &= \sum_{k=2}^s \sum_{u\in\mathcal{U}(k-2)}\lambda_u ([u_1,\dots,u_{k-2},x,y]-[u_1,\dots,u_{k-2},x,y]^\infty)_\ell\\
&= \sum_{k=2}^s \sum_{u\in\mathcal{U}(k-2)}\lambda_u \sum_{i_1+\dots i_k\le \ell-1}[(u_1)_{i_1},\dots,(u_{k-2})_{i_{k-2}},x_{i_{k-1}},y_{i_k}]_\ell.
\end{align*}

By definition of $\sigma$, we have $\|(u_j)_{i_j}\|\le M^{i_j}$ for all $j$. By the inequality $[u,v]\le\|u\|\|v\|$, for each term we deduce
$$\|[(u_1)_{i_1},\dots,(u_{k-2})_{i_{k-2}},x_{i_{k-1}},y_{i_k}]_\ell\|\le M^k\le M^{\ell-1},$$
so for some constant $C$ we have
$$\|(xy-x\boxtimes y)_\ell\|\le CM^{\ell-1}.\qedhere$$
\end{proof}

\begin{prop}\label{estim}
We have 
$$\sup_{x,y}\frac{\sigma((xy)^{-1}(x\boxtimes y))}{\max(\sigma(x),\sigma(y))^{1-1/s}}<\infty;$$
in particular, we have 
$$\frac{\sigma((xy)^{-1}(x\boxtimes y))}{\max(\sigma(x),\sigma(y))}\to 0\quad\text{when }\max(\sigma(x),\sigma(y))\to+\infty.$$
\end{prop}
\begin{proof}
Let us write the Baker-Campbell-Hausdorff formula as in (\ref{bch}).
\begin{align}\label{uuuu}&((xy)^{-1}(x\boxtimes y))_\ell= (-xy+x\boxtimes y)_\ell\\
&+\sum_{k=2}^s \sum_{u\in\mathcal{U}(k-2)[-xy,x\boxtimes y]}\lambda_u [u_1,\dots,u_{k-2},-xy,-xy+x\boxtimes y]_\ell.\notag\end{align}

Fix $u\in\mathcal{U}(k-2)[-xy,x\boxtimes y]$ and write $u_{k-1}=-xy$. Expanding, we obtain
\begin{align}\label{uuuuuu}
& [u_1,\dots,u_{k-1},-xy+x\boxtimes y]_\ell \\ & = \sum_{i_1+\dots +i_k= \ell}[(u_1)_{i_1},\dots,(u_{k-1})_{i_{k-1}},(-xy+x\boxtimes y)_{i_k}]_\ell.\notag
\end{align}
Since $\sigma(u_j)\le M$ for all $j<k$, we have $\|(u_j)_{i_j}\|\le M^{i_j}$. Using the inequality $\|[u,v]\|\le\|u\|v\|$ we deduce
$$\|[(u_1)_{i_1},\dots,(u_{k-1})_{i_{k-1}},(-xy+x\boxtimes y)_{i_k}]_\ell\|$$ $$\le M^{i_1+\dots+i_{k-1}}\|(-xy+x\boxtimes y)_{i_k}\|.$$
By Lemma \ref{compar}, we have $\|(-xy+x\boxtimes y)_{i}\|\le C_1M^{i-1}$ for some constant $C_1$ (chosen independent of $i$), so if $i_1+\dots+i_k=\ell$, using Lemma \ref{compar} we deduce
$$\|[(u_1)_{i_1},\dots,(u_{k-1})_{i_{k-1}},(-xy+x\boxtimes y)_{i_k}]_\ell \|\le C_1M^{\ell-1}.$$
Therefore from (\ref{uuuuuu}) we deduce, for some constant $C_2$
$$\|[u_1,\dots,u_{k-1},-xy+x\boxtimes y]_\ell\|\le C_2M^{\ell-1};$$
then from (\ref{uuuu}) we obtain for some other constant $C_3$
$$\|((xy)^{-1}(x\boxtimes y))_\ell\|\le C_3M^{\ell-1},$$
hence
$$\|((xy)^{-1}(x\boxtimes y))_\ell\|^{1/\ell}\le C_3^{1/\ell}M^{1-1/\ell}$$
so that, for some constant $C$
$$\sigma((xy)^{-1}(x\boxtimes y))\le C\max_{1\le\ell\le s} M^{1-1/\ell}=CM^{1-1/s}.\qedhere$$
\end{proof}


\baselineskip=16pt

\end{document}